\documentclass{amsart}

\usepackage{amssymb}
\usepackage{amsmath}
\usepackage{amsthm}
\usepackage{graphics}

\allowdisplaybreaks

\numberwithin{equation}{section}

\theoremstyle{plain}
\newtheorem{thm}{Theorem}[section]

\newtheorem{lem}[thm]{Lemma}

\theoremstyle{definition}

\newcommand{\R}{\mathbb{R}}

\newcommand{\Z}{\mathbb{Z}}

\newcommand{\calF}{\mathcal{F}}
\newcommand{\calL}{\mathcal{L}}
\newcommand{\calS}{\mathcal{S}}

\begin{document}

\title[Pseudo-differential operators on modulation spaces]
{A counterexample for boundedness \\ of pseudo-differential operators
\\ on modulation spaces}
\author{Mitsuru Sugimoto \and Naohito Tomita}
\date{}

\address{Mitsuru Sugimoto \\
Department of Mathematics \\
Graduate School of Science \\
Osaka University \\
Toyonaka, Osaka 560-0043, Japan}
\email{sugimoto@math.wani.osaka-u.ac.jp}

\address{Naohito Tomita \\
Department of Mathematics \\
Graduate School of Science \\
Osaka University \\
Toyonaka, Osaka 560-0043, Japan}
\email{tomita@gaia.math.wani.osaka-u.ac.jp}



\begin{abstract}
We prove that pseudo-differential operators with symbols in the class 
$S_{1,\delta}^0$ ($0<\delta<1$) are not always bounded on 
the modulation space $M^{p,q}$ ($q\neq2$).
\end{abstract}

\maketitle

\section{Introduction}\label{section1}

In 1980's, by Feichtinger
\cite{Feichtinger1,Feichtinger2},
the modulation spaces $M^{p,q}$ were introduced
as a fundamental function space of time-frequency
analysis which is originated in signal analysis or quantum mechanics.
See also \cite{Feichtinger3} or Triebel \cite{Triebel}.
The exact definition will be given in the next section, but
the main idea is to consider the decaying property
of a function with respect to the space variable and the
variable of its Fourier transform simultaneously.
Based on the same idea, Sj\"ostrand \cite{Sjostrand}
independently introduced a symbol class which assures the $L^2$-boundedness
of corresponding pseudo-differential operators.
Since this pioneering work,
the modulation spaces have been also recognized as a
useful tool for the theory of pseudo-differential operators
(see Gr\"ochenig \cite{Grochenig2}).
\par
In this paper, we investigate the boundedness properties of
pseudo-differential operators with symbols in $S_{\rho,\delta}^m$
on the modulation spaces $M^{p,q}$.
There have been already several literatures on this subject.
For example, Gr\"ochenig and Heil \cite{Grochenig-Heil},
Tachizawa \cite{Tachizawa},
Toft \cite{Toft}
proved that pseudo-differential operators
with symbols in $S_{0,0}^0$ are $M^{p,q}$-bounded.
On the other hand, Calder\'on and Vaillancourt \cite{C-V} proved that
pseudo-differential operators
with symbols in $S_{\delta,\delta}^0$ with $0<\delta<1$
(hence $S_{1,\delta}^0$) are $L^2$-bounded (hence $M^{2,2}$-bounded)
by reducing it to the case of $S_{0,0}^0$.
In view of these results, the class $S_{1,\delta}^0$ with $0<\delta<1$
appears to induce the $M^{p,q}$-boundedness, as well.
The objective of this paper is to show that this is not true:
\begin{thm}\label{1.1}
Let $1<p,q<\infty$, $m \in \R$ and $0<\delta<1$.
If $m>-|1/q-1/2|\delta n$,
then there exists a symbol $\sigma \in S_{1,\delta}^m$ such that
$\sigma(X,D)$ is not bounded on $M^{p,q}(\R^n)$.
\end{thm}
In particular, Theorem \ref{1.1} actually says that
symbols in the class $S_{1,\delta}^0$ ($0<\delta<1$)
do not always induce the $M^{p,q}$-boundedness in the case $q\neq2$.
We will prove this fact by constructing a counterexample.

\section{Main Result}\label{section2}
Let $\calS(\R^n)$ and $\calS'(\R^n)$ be the Schwartz spaces of
rapidly decreasing smooth functions
and tempered distributions,
respectively.
We define the Fourier transform $\calF f$
and the inverse Fourier transform $\calF^{-1}f$
of $f \in \calS(\R^n)$ by
\[
\calF f(\xi)
=\widehat{f}(\xi)
=\int_{\R^n}e^{-i\xi \cdot x}\, f(x)\, dx
\quad \text{and} \quad
\calF^{-1}f(x)
=\frac{1}{(2\pi)^n}
\int_{\R^n}e^{ix\cdot \xi}\, f(\xi)\, d\xi.
\]

Let $m \in \R$ and $0 \le \delta \le \rho \le 1$.
The symbol class $S_{\rho,\delta}^m$ consists of
all $\sigma \in C^{\infty}(\R^n \times \R^n)$ such that
\[
|\partial_{\xi}^{\alpha}\partial_x^{\beta}\sigma(x,\xi)|
\le C_{\alpha,\beta}(1+|\xi|)^{m-\rho |\alpha|+\delta |\beta|}
\]
for all $\alpha,\beta \in \Z_+^n=\{0,1,2,\dots\}^n$.
For $\sigma \in S_{\rho,\delta}^m$,
we define the pseudo-differential operator $\sigma(X,D)$
by
\[
\sigma(X,D)f(x)
=\frac{1}{(2\pi)^n}
\int_{\R^n}e^{ix\cdot \xi}\,
\sigma(x,\xi)\, \widehat{f}(\xi)\, d\xi
\]
for $f \in \calS(\R^n)$.
Given a symbol $\sigma \in S_{\rho,\delta}^m$
with $\delta<1$,
the symbol $\sigma^*$ defined by
\begin{align}\label{(2.1)}
\sigma^*(x,\xi)
&=\textrm{Os-}
\frac{1}{(2\pi)^n}\int_{\R^n}\int_{\R^n}
e^{-iy\cdot \zeta}\,
\overline{\sigma(x+y,\xi+\zeta)}\,
dy \, d\zeta \\
&=\lim_{\epsilon \to 0}
\frac{1}{(2\pi)^n}
\int_{\R^n}\int_{\R^n}
e^{-iy\cdot \zeta}\,
\chi(\epsilon y, \epsilon \zeta)\,
\overline{\sigma(x+y,\xi+\zeta)}\,
dy \, d\zeta \nonumber
\end{align}
satisfies $\sigma^* \in S_{\rho,\delta}^m$
and
\begin{equation}\label{(2.2)}
(\sigma(X,D)f,g)=(f,\sigma^*(X,D)g)
\qquad \text{for all} \ f,g \in \calS(\R^n),
\end{equation}
where $\chi \in \calS(\R^{2n})$
satisfies $\chi(0,0)=1$
and $(\cdot, \cdot )$
denotes the inner product on $L^2(\R^n)$
(\cite[Chapter 2, Theorem 2.6]{Kumano-go}).
Note that oscillatory integrals
are independent of the choice of $\chi \in \calS(\R^{2n})$
satisfying $\chi(0,0)=1$
(\cite[Chapter 1, Theorem 6.4]{Kumano-go}).

We introduce the modulation spaces
based on Gr\"ochenig \cite{Grochenig}.
Fix a function $\gamma \in \calS(\R^n)\setminus \{ 0 \}$
(called the window function).
Then the short-time Fourier transform $V_{\gamma}f$ of
$f \in \calS'(\R^n)$ with respect to $\gamma$
is defined by
\[
V_{\gamma}f(x,\xi)
=(f, M_{\xi}T_x \gamma)
\qquad \text{for} \ x, \xi \in \R^n,
\]
where $M_{\xi}\gamma(t)=e^{i\xi \cdot t}\gamma(t)$
and $T_x \gamma(t)=\gamma(t-x)$.
We can express it in a form of the integral
\[
V_{\gamma}f(x,\xi)
=\int_{\R^n}f(t)\, \overline{\gamma(t-x)}\,
e^{-i\xi \cdot t}\, dt,
\]
which has actually the meaning
for an appropriate function $f$ on $\R^n$.
We note that,
for $f \in \calS'(\R^n)$,
$V_{\gamma}f$ is continuous on $\R^{2n}$
and $|V_{\gamma}f(x,\xi)|\le C(1+|x|+|\xi|)^N$
for some constants $C,N \ge 0$
(\cite[Theorem 11.2.3]{Grochenig}).
Let $1\le p,q \le \infty$.
Then the modulation space $M^{p,q}(\R^n)$
consists of all $f \in \calS'(\R^n)$
such that
\[
\|f\|_{M^{p,q}}
=\|V_{\gamma}f\|_{L^{p,q}}
=\left\{ \int_{\R^n} \left(
\int_{\R^n} |V_{\gamma}f(x,\xi)|^{p}\, dx
\right)^{q/p} d\xi \right\}^{1/q}
< \infty.
\]
We note that
$M^{2,2}(\R^n)=L^2(\R^n)$
(\cite[Proposition 11.3.1]{Grochenig})
and $M^{p,q}(\R^n)$ is a Banach space
(\cite[Proposition 11.3.5]{Grochenig}).
The definition of $M^{p,q}(\R^n)$ is independent
of the choice of the window function
$\gamma \in \calS(\R^n)\setminus \{ 0 \}$,
that is,
different window functions
yield equivalent norms
(\cite[Proposition 11.3.2]{Grochenig}).
\par
We also introduce a special symbol which will act as the counterexample
for the boundedness stated in Introduction.
Let $\varphi,\eta \in \calS(\R^n)$ be real-valued functions
satisfying
\begin{align*}
\varphi:
\quad &{\rm supp}\, \varphi \subset \{\xi: |\xi| \le 1/8\},
\quad \int_{\R^n}\varphi(\xi)\, d\xi =1, \\
\eta:
\quad &{\rm supp}\, \eta \subset
\{\xi: 2^{-1/2}\le |\xi| \le 2^{1/2}\},
\quad \eta=1 \ \text{on} \
\{\xi: 2^{-1/4}\le |\xi| \le 2^{1/4}\}.
\end{align*}
Moreover, we assume that $\varphi$ is radial.
Then we define the symbol $\sigma_\delta$ by
\begin{equation}\label{(2.3)}
\sigma_{\delta}(x,\xi)
=\sum_{j=j_0}^{\infty}2^{jm}
\left(\sum_{0<|k|\le 2^{j\delta/2}}
e^{-ik\cdot(2^{j\delta/2}x-k)}\,
\Phi(2^{j\delta/2}x-k) \right)
\eta(2^{-j}\xi),
\end{equation}
where $\Phi=\calF^{-1}\varphi$,
$0<\delta<1$ and $j_0\in\Z_+$ is chosen to satisfy
\[
1+2^{j_0(\delta-1)+1} \le 2^{1/4}, \quad
1-2^{j_0(\delta-1)+1} \ge 2^{-1/4}, \quad
2^{-j_0 \delta/2} \sqrt{n} \le 2^{-3}.
\]
The symbol $\sigma_\delta^*$ is constructed from
$\sigma_\delta$ using the oscillatory integral \eqref{(2.1)}.

\medskip
Now, we state our main result which is a precise version of
Theorem \ref{1.1} in Introduction.
\begin{thm}\label{2.1}
Let $1<p,q<\infty$, $0<\delta<1$, and $m>-|1/q-1/2|\delta n$.
Then the symbols $\sigma_\delta$ and $\sigma_\delta^*$
defined by \eqref{(2.3)}
belong to the class $S_{1,\delta}^m$.
Moreover,
if $q\geq2$ {\rm (}$q\leq2$ resp.{\rm )}, then the corresponding operator
$\sigma_\delta(X,D)$ {\rm (}$\sigma_\delta^*(X,D)$ resp.{\rm )} is
not bounded on $M^{p,q}(\R^n)$.
\end{thm}
The proof of Theorem \ref{2.1} will be given in the next section.

\section{Proof}\label{section3}
In the below, we consider the symbol $\tau_{\delta}$
instead of $\sigma_\delta$ for the sake of simplicity.
In order to avoid confusion, we repeat the notation in this context,
and also introduce a family of functions $\{f_{j,\epsilon,\delta}\}_j$.
Let $\varphi,\psi,\eta \in \calS(\R^n)$ be
real-valued functions satisfying
\begin{align*}
\varphi:
\quad &{\rm supp}\, \varphi \subset \{\xi: |\xi| \le 1/8\},
\quad \int_{\R^n}\varphi(\xi)\, d\xi =1, \\
\psi:
\quad &{\rm supp}\, \psi \subset \{\xi: |\xi| \le 1/2\},
\quad \psi=1 \ \text{on} \ \{\xi: |\xi| \le 1/4\}, \\
\eta:
\quad &{\rm supp}\, \eta \subset
\{\xi: 2^{-1/2}\le |\xi| \le 2^{1/2}\},
\quad \eta=1 \ \text{on} \
\{\xi: 2^{-1/4}\le |\xi| \le 2^{1/4}\}.
\end{align*}
Moreover, we assume that
$\varphi$ and $\psi$ are radial.
This assumption implies that
$\Phi$ and $\Psi$ are also real-valued functions,
where $\Phi=\calF^{-1}\varphi$ and $\Psi=\calF^{-1}\psi$.
For these $\Phi, \Psi, \eta$,
we define the symbol $\tau_\delta$
and the functions $f_{j,\epsilon,\delta}$ by
\begin{equation}\label{(3.1)}
\tau_\delta(x,\xi)
=\sum_{j=j_0}^{\infty}2^{jm}\left(
\sum_{0<|k|\le 2^{j\delta}}
e^{-ik\cdot(2^{j\delta}x-k)}\,
\Phi(2^{j\delta}x-k) \right)
\eta(2^{-j}\xi)
\end{equation}
and
\begin{equation}\label{(3.2)}
f_{j,\epsilon,\delta}(x)=
\sum_{0<|k'|\le 2^{j\delta}}
|k'|^{-n/q-\epsilon} \,
e^{ik'\cdot(x-k')}\,
\Psi(x-k'),
\end{equation}
where $k,k' \in \Z^n$,
$\epsilon>0$,
$0<\delta<1/2$ and $j_0 \in \Z_+$ is chosen to satisfy
\[
1+2^{j_0(2\delta-1)+1} \le 2^{1/4}, \quad
1-2^{j_0(2\delta-1)+1} \ge 2^{-1/4}, \quad
2^{-j_0 \delta} \sqrt{n} \le 2^{-3}.
\]
Note that $\sigma_{2\delta}(x,\xi)=\tau_\delta(x,\xi)$.
\begin{lem}\label{3.1}
The symbol $\tau_\delta$ defined by {\rm \eqref{(3.1)}}
belongs to $S_{1,2\delta}^m$.
\end{lem}
\begin{proof}
Since
${\rm supp}\, \eta(2^{-j}\cdot)
\subset \{2^{j-1/2} \le |\xi| \le 2^{j+1/2}\}$,
we see that,
for each $\xi \in \R^n$,
at most one term in the sum \eqref{(3.1)} is nonzero
with respect to $j$.
Note that $2^j \sim |\xi| \sim 1+|\xi|$
on ${\rm supp}\, \eta(2^{-j}\cdot)$.
Let $\alpha, \beta \in \Z_+^n$
and $\xi \in {\rm supp}\, \eta(2^{-j}\cdot)$.
Using
\begin{multline*}
\partial_{\xi}^{\alpha}\partial_x^{\beta}
\tau_\delta(x,\xi)
=2^{jm}
\sum_{0<|k|\le 2^{j\delta}}
\sum_{\beta_1+\beta_2=\beta}
C_{\beta_1,\beta_2}\,
(2^{j\delta}k)^{\beta_1}\,
e^{-ik\cdot(2^{j\delta}x-k)} \\
\times  2^{j\delta|\beta_2|}\,
(\partial^{\beta_2}\Phi)(2^{j\delta}x-k)\,
2^{-j|\alpha|}
(\partial^{\alpha}\eta)(2^{-j}\xi),
\end{multline*}
we have
\begin{align*}
|\partial_{\xi}^{\alpha}\partial_x^{\beta} \tau_\delta(x,\xi)|
&\le C2^{j(m-|\alpha|)}
\| \partial^{\alpha}\eta \|_{L^{\infty}}
\left( \sum_{0<|k|\le 2^{j\delta}}(1+|2^{j\delta}x-k|)^{-n-1} \right) \\
&\qquad \qquad \times
\left( \sum_{\beta_1+\beta_2=\beta}
2^{j\delta(2|\beta_1|+|\beta_2|)}
\| (1+|\cdot|)^{n+1}(\partial^{\beta_2}\Phi)\|_{L^{\infty}}
\right) \\
&\le C2^{j(m-|\alpha|+2\delta|\beta|)}
\left( \sup_{y \in \R^n}
\sum_{k \in \Z^n}(1+|y-k|)^{-n-1} \right)  \\
&\le C(1+|\xi|)^{m-|\alpha|+2\delta|\beta|}.
\end{align*}
In the case
$\xi \not\in \cup_{j \ge j_0} \, {\rm supp}\, \eta(2^{-j}\cdot)$,
we have nothing to prove.
\end{proof}
\begin{lem}\label{3.2}
Let $1<p,q<\infty$,
and $f_{j,\epsilon,\delta}$ be defined by {\rm \eqref{(3.2)}}.
Then the following are true:
\begin{enumerate}
\item
The Fourier transform of
$e^{i2^{j(1-\delta)}x_1}\, f_{j,\epsilon,\delta}(x)$ is
\[
\calF[M_{2^{j(1-\delta)}e_1}f_{j,\epsilon,\delta}](\xi)
=\sum_{0<|k'|\le 2^{j\delta}}
|k'|^{-n/q-\epsilon}\,
e^{-ik' \cdot (\xi-2^{j(1-\delta)}e_1)}\,
\psi(\xi-2^{j(1-\delta)}e_1-k'),
\]
where $x=(x_1,x_2, \dots, x_n) \in \R^n$
and $e_1=(1,0,\dots,0) \in \R^n$.
\item
There exists a constant $C>0$ such that
\[
\|f_{j,\epsilon,\delta}(2^{j\delta}\cdot)\|_{M^{p,q}}
\le C2^{j\delta n(1/q-1)}
\qquad \text{for all} \  j \ge j_0.
\]
\end{enumerate}
\end{lem}
\begin{proof}
We consider only (2).
Let $\gamma \in \calS(\R^n)\setminus\{0\}$.
Since
\begin{align*}
&V_{\gamma}[f_{j,\epsilon,\delta}(2^{j\delta}\cdot)]
(2^{-j\delta}x,2^{j\delta}\xi) \\
&=\sum_{0<|k'|\le 2^{j\delta}}
|k'|^{-n/q-\epsilon} \,
e^{-i|k'|^2}\, \int_{\R^n}e^{ik'\cdot(2^{j\delta}t)}\,
\Psi(2^{j\delta}t-k')\,
\overline{\gamma(t-2^{-j\delta}x)}\,
e^{-i(2^{j\delta}\xi)\cdot t}\, dt \\
&=2^{-j\delta n}
\sum_{0<|k'|\le 2^{j\delta}}
|k'|^{-n/q-\epsilon} \,
e^{-i|k'|^2}\, \int_{\R^n}e^{-i(\xi-k')\cdot t}\,
\Psi(t-k')\,
\overline{\gamma(2^{-j\delta}(t-x))}\, dt \\
&=2^{-j\delta n}
\sum_{0<|k'|\le 2^{j\delta}}
|k'|^{-n/q-\epsilon} \, e^{-i|k'|^2} \\
&\qquad \times
\int_{\R^n}
\left\{ (1+|\xi-k'|^2)^{-n}\, (I-\Delta_t)^n\,
e^{-i(\xi-k')\cdot t}\right\}
\Psi(t-k')\,
\overline{\gamma(2^{-j\delta}(t-x))}\, dt,
\end{align*}
we have
\begin{multline*}
|V_{\gamma}[f_{j,\epsilon,\delta}(2^{j\delta}\cdot)]
(2^{-j\delta}x,2^{j\delta}\xi)|
\le C2^{-j\delta n}
\sum_{|\alpha_1+\alpha_2| \le 2n}
\sum_{0<|k'|\le 2^{j\delta}}
|k'|^{-n/q-\epsilon} \\
\times (1+|\xi-k'|)^{-2n}
[|\widetilde{\partial^{\alpha_1}\Psi}|*
|(\partial^{\alpha_2}\gamma)(2^{-j\delta}\cdot)|](k'-x),
\end{multline*}
where $\widetilde{\Psi}(t)=\Psi(-t)$.
Hence,
we get
\begin{align*}
&\|V_{\gamma}[f_{j,\epsilon,\delta}
(2^{j\delta}\cdot)]\|_{L^{p,q}} \\
&=2^{-j\delta n(1/p-1/q)}\left\{ \int_{\R^n}
\left( \int_{\R^n}
\left| V_{\gamma}[f_{j,\epsilon,\delta}(2^{j\delta}\cdot)]
(2^{-j\delta}x,2^{j\delta}\xi)
\right|^p dx \right)^{q/p} d\xi \right\}^{1/q} \\
&\le C2^{-j\delta n(1/p-1/q+1)}
\sum_{|\alpha_1+\alpha_2| \le 2n} \Bigg\{ \int_{\R^n} \\
&\quad \times \Bigg(
\sum_{0<|k'|\le 2^{j\delta}}
|k'|^{-n/q-\epsilon}\, (1+|\xi-k'|)^{-2n}
\| |\widetilde{\partial^{\alpha_1}\Psi}|*
|(\partial^{\alpha_2}\gamma)(2^{-j\delta}\cdot)|
\|_{L^p} \Bigg)^q d\xi \Bigg\}^{1/q} \\
&\le C2^{-j\delta n(1/p-1/q+1)}
\Bigg(\sum_{|\alpha_1+\alpha_2| \le 2n}
\|\widetilde{\partial^{\alpha_1}\Psi}\|_{L^1}
\|(\partial^{\alpha_2}\gamma)(2^{-j\delta}\cdot)\|_{L^p}\Bigg) \\
&\quad \times
\Bigg\{ \sum_{k \in \Z^n}
\int_{k+[-1/2,1/2]^n}
\Bigg( \sum_{k' \neq 0}
|k'|^{-n/q-\epsilon}\, (1+|\xi-k'|)^{-2n} \Bigg)^q
d\xi \Bigg\}^{1/q} \\
&\le C2^{j\delta n(1/q-1)}
\Bigg\{ \sum_{k \in \Z^n}
\Bigg( \sum_{k' \neq 0}
|k'|^{-n/q-\epsilon}\,
(1+|k-k'|)^{-2n} \Bigg)^q \Bigg\}^{1/q}
\le C2^{j\delta n(1/q-1)}.
\end{align*}
The proof is complete.
\end{proof}
\begin{lem}\label{3.3}
Let $\tau_\delta$ be defined by {\rm \eqref{(3.1)}},
and
\begin{equation}\label{(3.3)}
g_{j,\epsilon,\delta}(x)
=\sum_{0<|k|,|k'|\le 2^{j\delta}} 
|k'|^{-n/q-\epsilon} \,
e^{-ik\cdot(x-k)}\, e^{ik' \cdot (x-k')}\,
\Phi(x-k)\, \Psi(x- k').
\end{equation}
Then
\[
\tau_\delta(X,D)
[(M_{2^{j(1-\delta)}e_1}f_{j,\epsilon,\delta})(2^{j\delta}\cdot)](x)
=2^{jm}e^{i 2^j x_1}g_{j,\epsilon,\delta}(2^{j\delta}x)
\]
for all $j \ge j_0$ and $x \in \R^n$,
where $f_{j,\epsilon,\delta}$ is defined by {\rm \eqref{(3.2)}}.
\end{lem}
\begin{proof}
By Lemma \ref{3.2} (1),
we have
\begin{align*}
&\calF[(M_{2^{j(1-\delta)}e_1}
f_{j,\epsilon,\delta})(2^{j\delta}\cdot)](\xi)
=2^{-j\delta n}
\calF[M_{2^{j(1-\delta)}e_1}f_{j,\epsilon,\delta}]
(2^{-j\delta}\xi) \\
&=\sum_{0<|k'|\le 2^{j\delta}}
|k'|^{-n/q-\epsilon} \,
e^{-ik' \cdot (2^{-j\delta}\xi-2^{j(1-\delta)}e_1)}\,
2^{-j\delta n}\psi(2^{-j\delta}\xi-(2^{j(1-\delta)}e_1+k')).
\end{align*}
Since $2\delta<1$,
$j \ge j_0$,
$1+2^{j_0(2\delta-1)+1} \le 2^{1/4}$
and $1-2^{j_0(2\delta-1)+1} \ge 2^{-1/4}$,
we see that
\begin{align*}
{\rm supp}\, \psi(2^{-j\delta}\cdot-(2^{j(1-\delta)}e_1+k'))
&\subset
\{\xi: |\xi-(2^j e_1+2^{j\delta}k')| \le 2^{j\delta-1}\} \\
&\subset
\{\xi: 2^j-2^{2j\delta+1}\le |\xi| \le 2^j+2^{2j\delta+1}\} \\
&\subset \{\xi: 2^{j-1/4}\le |\xi| \le 2^{j+1/4}\}
\end{align*}
for all $|k'| \le 2^{j\delta}$.
This implies
\[
\mathrm{supp}\, \calF[(M_{2^{j(1-\delta)}e_1}
f_{j,\epsilon,\delta})(2^{j\delta}\cdot)]
\subset \{\xi: 2^{j-1/4} \le \xi \le 2^{j+1/4}\}
\]
for all $j \ge j_0$.
Hence, noting that
${\rm supp}\, \eta(2^{-j'}\cdot)
\subset \{2^{j'-1/2}\le |\xi| \le 2^{j'+1/2}\}$
and $\eta(2^{-j'}\cdot)=1$
on $\{2^{j'-1/4}\le |\xi| \le 2^{j'+1/4}\}$,
we obtain
\begin{align*}
&\tau_\delta(X,D)
[(M_{2^{j(1-\delta)}e_1}f_{j,\epsilon,\delta})(2^{j\delta}\cdot)](x) \\
&=\frac{1}{(2\pi)^n}
\int_{\R^n}e^{ix\cdot\xi}
\left\{2^{jm}\left(
\sum_{0<|k|\le 2^{j\delta}}
e^{-ik\cdot(2^{j\delta}x-k)}\,
\Phi(2^{j\delta}x-k) \right)
\eta(2^{-j}\xi)\right\} \\
&\qquad \times
\calF[(M_{2^{j(1-\delta)}e_1}
f_{j,\epsilon,\delta})(2^{j\delta}\cdot)](\xi)\, d\xi \\
&=2^{jm}\left(
\sum_{0<|k|\le 2^{j\delta}}
e^{-ik\cdot(2^{j\delta}x-k)}\,
\Phi(2^{j\delta}x-k) \right) \\
&\qquad \times \frac{1}{(2\pi)^n}
\int_{\R^n}e^{ix\cdot\xi}\,
\calF[(M_{2^{j(1-\delta)}e_1}
f_{j,\epsilon,\delta})(2^{j\delta}\cdot)](\xi)\, d\xi \\
&=2^{jm}\left(
\sum_{0<|k|\le 2^{j\delta}}
e^{-ik\cdot(2^{j\delta}x-k)}\,
\Phi(2^{j\delta}x-k) \right)
(M_{2^{j(1-\delta)}e_1}
f_{j,\epsilon,\delta})(2^{j\delta}x) \\
&=2^{jm}e^{i2^j x_1}\, g_{j,\epsilon,\delta}(2^{j\delta}x)
\end{align*}
for all $j \ge j_0$ and $x \in \R^n$.
The proof is complete.
\end{proof}
\begin{lem}[{\cite[Corollary 11.2.7]{Grochenig}}]\label{3.4}
Let $f \in \calS'(\R^n)$ and $\gamma \in \calS(\R^n)\setminus\{0\}$.
Then
\[
(f,g)
=\frac{1}{\|\gamma\|_{L^2}^2}
\int_{\R^{2n}}
V_{\gamma}f(x,\xi)\, \overline{V_{\gamma}g(x,\xi)}\,
dx \, d\xi
\qquad
\text{for all} \ g \in \calS(\R^n).
\]
\end{lem}
For $1 \le p \le \infty$,
$p'$ is the conjugate exponent of $p$
(that is, $1/p+1/p'=1$).
\begin{lem}[{\cite[Proposition 11.3.4, Theorem 11.3.6]{Grochenig}}]
\label{3.5}
If $1 \le p,q <\infty$,
then $\calS(\R^n)$ is dense in $M^{p,q}(\R^n)$
and $M^{p,q}(\R^n)^*=M^{p',q'}(\R^n)$
under the duality
\[
\langle f,g\rangle
=\frac{1}{\|\gamma\|_{L^2}^2}
\int_{\R^{2n}}
V_{\gamma}f(x,\xi)\, \overline{V_{\gamma}g(x,\xi)}\,
dx \, d\xi
\]
for $f \in M^{p,q}(\R^n)$ and $g \in M^{p',q'}(\R^n)$,
where $\gamma \in \calS(\R^n)\setminus\{0\}$.
\end{lem}
We denote by $B$ the tensor product of B-spline of degree $2$,
that is,
\[
B(t)=\prod_{i=1}^n \chi_{[-1/2,1/2]}*\chi_{[-1/2,1/2]}(t_i),
\]
where
$\chi_{[-1/2, 1/2]}$ is the characteristic function of $[-1/2,1/2]$.
Note that
${\rm supp}\, B \subset [-1,1]^n$
and
$\calF^{-1} B(t)
=(2\pi)^{-n}\prod_{i=1}^n\{(\sin (t_i/2))/(t_i/2)\}^2 \in M^{p,q}(\R^n)$
for all $1 \le p,q \le \infty$.
By Lemmas \ref{3.4} and \ref{3.5},
if $1<p,q< \infty$ and $f \in \calS(\R^n)$,
then
\begin{align}\label{(3.4)}
&\|f\|_{M^{p,q}}
=\sup_{\|g\|_{M^{p',q'}}=1}
\left| \langle f,g \rangle \right| \ge
\left| \left\langle f,  \frac{\calF^{-1}B}{\|\calF^{-1}B\|_{M^{p',q'}}}
\right\rangle \right| \\
&=\frac{1}{\|\calF^{-1}B\|_{M^{p',q'}}}
\left|\int_{\R^{2n}}
V_{\gamma}f(x,\xi)\, \overline{V_{\gamma}[\calF^{-1}B](x,\xi)}
\, dx \, d\xi \right| \nonumber \\
&=\frac{1}{\|\calF^{-1}B\|_{M^{p',q'}}}
\left| \int_{\R^n}f(t)\, \calF^{-1}B(t)\, dt\right|, \nonumber
\end{align}
where $\gamma \in \calS(\R^n)$ such that
$\|\gamma\|_{L^2}=1$.
\begin{lem}\label{3.6}
Let $1<p,q<\infty$,
and $g_{j,\epsilon,\delta}$ be defined by {\rm \eqref{(3.3)}}.
Then there exists a constant $C>0$ such that
\[
\|g_{j,\epsilon,\delta} (2^{j\delta}\cdot)\|_{M^{p,q}}
\ge C2^{-j\delta(n/q+\epsilon)}
\qquad \text{for all} \ j \ge j_0.
\]
\end{lem}
\begin{proof}
Let $B$ be the tensor product of B-spline of degree $2$.
Note that
$g_{j,\epsilon,\delta}(2^{j\delta}\cdot) \in \calS(\R^n)$.
By \eqref{(3.4)},
we have
\begin{align*}
&\|g_{j,\epsilon,\delta} (2^{j\delta}\cdot)\|_{M^{p,q}}
\ge \frac{1}{\|\calF^{-1}B\|_{M^{p',q'}}}
\left| \int_{\R^n}
g_{j,\epsilon,\delta}(2^{j\delta}x)\,
\calF^{-1}B(x)\, dx \right| \\
&=\frac{1}{\|\calF^{-1}B\|_{M^{p',q'}}}
\bigg| \bigg\{ \sum_{0<|k|\le 2^{j\delta}}
 |k|^{-n/q-\epsilon} \,
\int_{\R^n}\Phi(2^{j\delta}x-k)\, \Psi(2^{j\delta}x-k)\,
\calF^{-1}B(x)\, dx  \bigg\} \\
&\quad + \bigg\{ \sum_{0<|k|\le 2^{j\delta}}
\sum_{\scriptstyle 0<|k'| \le 2^{j\delta}
\atop \scriptstyle k' \neq k}
|k'|^{-n/q-\epsilon} \\
&\quad \quad \times
\int_{\R^n}\Big( e^{-ik\cdot (2^{j\delta}x-k)}\,
\Phi(2^{j\delta}x-k)\, 
\calF^{-1}B(x)\Big)
\Big( e^{ik'\cdot(2^{j\delta}x-k')}\,
\Psi(2^{j\delta}x- k')\Big) \, dx \bigg\}\bigg| \\
&=\frac{1}{\|\calF^{-1}B\|_{M^{p',q'}}}|I+II|.
\end{align*}
We first consider $I$.
Note that our assumptions $\int_{\R^n}\varphi(\xi)\, d\xi=1$
, ${\rm supp}\, \varphi \subset \{|\xi| \le 1/8\}$
and $\psi=1$ on $\{|\xi| \le 1/4\}$
give $\varphi*\psi=1$ on $\{|\xi| \le 1/8\}$.
Since
${\rm supp}\, B \subset \{ |\xi| \le \sqrt{n}\}$ and
$2^{-j\delta}\sqrt{n} \le 2^{-j_0 \delta}\sqrt{n} \le 1/8$,
by Plancherel's theorem,
we see that
\begin{align*}
&\sum_{0<|k|\le 2^{j\delta}}
|k|^{-n/q-\epsilon}\,
\int_{\R^n}\Phi(2^{j\delta}x-k)\, \Psi(2^{j\delta}x-k)\,
\calF^{-1}B(x)\, dx \\
&=\sum_{0<|k|\le 2^{j\delta}}
|k|^{-n/q-\epsilon}\,
\int_{\R^n}(\Phi\, \Psi)(2^{j\delta}x-k)\,
\overline{\calF^{-1}B(x)}\, dx \\
&=\sum_{0<|k|\le 2^{j\delta}}
|k|^{-n/q-\epsilon}\,
\frac{1}{(2\pi)^{2n}}
\int_{|\xi| \le \sqrt{n}}
2^{-j\delta n}\,
e^{-ik\cdot(2^{-j\delta}\xi)}\,
\varphi*\psi(2^{-j\delta}\xi)\, \overline{B(\xi)} \, d\xi \\
&=C_n2^{-j\delta n}\sum_{0<|k|\le 2^{j\delta}}
|k|^{-n/q-\epsilon}\,
\int_{\R^n}e^{-i(2^{-j\delta}k)\cdot\xi}\,
B(\xi) \, d\xi \\
&=C_n2^{-j\delta n}
\sum_{0<|k|\le 2^{j\delta}}
|k|^{-n/q-\epsilon}\,
\prod_{i=1}^n
\left( \frac{\sin k_i/2^{j\delta+1}}{k_i/2^{j\delta+1}} \right)^2.
\end{align*}
We next consider $II$.
Using
\[
\calF[T_k M_{-k}\Phi]=
e^{i|k|^2}T_{-k}M_{-k}\varphi
\quad \text{and} \quad
\calF[T_{k'}M_{-k'}\Psi]=
e^{i|k'|^2}T_{-k'}M_{-k'}\psi,
\]
we have
\begin{align*}
&\int_{\R^n}\left( e^{-ik\cdot(2^{j\delta}x-k)}\,
\Phi(2^{j\delta}x-k)\, 
\calF^{-1}B(x) \right)
\left( e^{ik'\cdot(2^{j\delta}x-k')}\,
\Psi(2^{j\delta}x- k')\right) \, dx \\
&=\int_{\R^n}
\left( T_k M_{-k}\Phi(x)\, 
2^{-j\delta n}[\calF^{-1}B](2^{-j\delta}x)\right)
\overline{T_{k'}M_{-k'}\Psi(x)} \, dx \\
&=\frac{e^{i(|k|^2-|k'|^2)}}{(2\pi)^{2n}}
\int_{\R^n}[(T_{-k}M_{-k}\varphi)*(B(2^{j\delta}\cdot))](\xi)\,
\overline{[T_{-k'}M_{-k'}\psi](\xi)}\, d\xi \\
&=\frac{e^{i(|k|^2-|k'|^2)}}{(2\pi)^{2n}}
\int_{\R^n}[(M_{-k}\varphi)*(B(2^{j\delta}\cdot))](\xi+k)\,
\overline{[M_{-k'}\psi](\xi+k')}\, d\xi.
\end{align*}
Since ${\rm supp}\, \varphi \subset \{|\xi| \le 1/8\}$,
${\rm supp}\, B(2^{j\delta}\cdot)
\subset \{|\xi| \le 2^{-j\delta}\sqrt{n}\}
\subset \{|\xi| \le 1/8\}$,
we see that
${\rm supp}\, (M_{-k}\varphi)*(B(2^{j\delta}\cdot))
\subset \{|\xi| \le 1/4\}$.
On the other hand,
${\rm supp}\, M_{-k'}\psi \subset \{|\xi|\le 1/2\}$.
Hence, if $k\neq k'$, then
\[
\frac{e^{i(|k|^2-|k'|^2)}}{(2\pi)^{2n}}
\int_{\R^n}[(M_{-k}\varphi)*(B(2^{j\delta}\cdot))](\xi+k)\,
\overline{[M_{-k'}\psi](\xi+k')}\, d\xi=0,
\]
that is, $II=0$.
Therefore,
since
$\prod_{i=1}^n(\sin x_i/x_i)^2 \ge C$ on $[-1/2,1/2]^n$,
we get
\begin{align*}
\frac{1}{\|\calF^{-1}B\|_{M^{p',q'}}}|I+II|
&=\frac{C_n 2^{-j\delta n}}{\|\calF^{-1}B\|_{M^{p',q'}}}
\left|\sum_{0<|k|\le 2^{j\delta}}
|k|^{-n/q-\epsilon}\,
\prod_{i=1}^n
\left( \frac{\sin k_i/2^{j\delta+1}}{k_i/2^{j\delta+1}} \right)^2 \right| \\
&=C_n 2^{-j\delta n}
\sum_{0<|k|\le 2^{j\delta}}
|k|^{-n/q-\epsilon}\,
\prod_{i=1}^n
\left( \frac{\sin k_i/2^{j\delta+1}}{k_i/2^{j\delta+1}} \right)^2 \\
&\ge C_n 2^{-j\delta n}\, 2^{-j\delta(n/q+\epsilon)}
\sum_{0<|k|\le 2^{j\delta}}1
\ge C2^{-j\delta(n/q+\epsilon)}.
\end{align*}
The proof is complete.
\end{proof}

We are now ready to prove Theorem \ref{2.1}.

\medskip
\noindent
{\it Proof of Theorem \ref{2.1}}.
Assume that $1<p,q<\infty$,
$0<\delta<1$ and $m>-|1/q-1/2|\delta n$.
Let $\sigma_\delta$ be defined by \eqref{(2.3)}.
Then $\sigma_\delta(x,\xi)=\tau_{\delta/2}(x,\xi)$,
where $\tau_{\delta/2}$ is defined by \eqref{3.1}.
By Lemma \ref{3.1},
we see that $\sigma_\delta \in S_{1,\delta}^m$.
This implies $\sigma_\delta^* \in S_{1,\delta}^m$,
where $\sigma_\delta^*$ is defined by \eqref{(2.1)}.

We first consider the case $q \ge 2$.
In this case, $m>(1/q-1/2)\delta n$.
Since $m>(1/q-1/2)\delta n$,
we can take $\epsilon>0$
such that $m-\epsilon \delta/2>(1/q-1/2)\delta n$.
We assume that $\sigma_\delta(X,D)$ is bounded on $M^{p,q}(\R^n)$.
Then, by Lemmas \ref{3.2}, \ref{3.3}, \ref{3.6}
and the modulation invariance of the norm $\| \cdot \|_{M^{p,q}}$,
we see that
\begin{align*}
C2^{j(m-\delta(n/q+\epsilon)/2)}
&\le 2^{jm}\|g_{j,\epsilon,\delta/2}(2^{j\delta/2}\cdot)\|_{M^{p,q}}
=\|2^{jm}e^{i 2^j x_1}
g_{j,\epsilon,\delta/2}(2^{j\delta/2}\cdot)\|_{M^{p,q}} \\
&=\|\tau_{\delta/2}(X,D)
[(M_{2^{j(1-\delta/2)}e_1}f_{j,\epsilon,\delta/2})(2^{j\delta/2}\cdot)]
\|_{M^{p,q}} \\
&=\|\sigma_{\delta}(X,D)
[(M_{2^{j(1-\delta/2)}e_1}f_{j,\epsilon,\delta/2})(2^{j\delta/2}\cdot)]
\|_{M^{p,q}} \\
&\le \|\sigma_{\delta}(X,D)\|_{\calL(M^{p,q})}
\|(M_{2^{j(1-\delta/2)}e_1}f_{j,\epsilon,\delta/2})
(2^{j\delta/2}\cdot)\|_{M^{p,q}} \\
&=\|\sigma_\delta(X,D)\|_{\calL(M^{p,q})}
\|f_{j,\epsilon,\delta/2}(2^{j\delta/2}\cdot)\|_{M^{p,q}}
\le C 2^{j\delta n(1/q-1)/2}
\end{align*}
for all $j \ge j_0$, where
$f_{j,\epsilon,\delta/2}$ and $g_{j,\epsilon,\delta/2}$
are defined by \eqref{(3.2)} and \eqref{(3.3)}.
However, since $m-\epsilon \delta/2>(1/q-1/2)\delta n$,
this is a contradiction.
Hence,
$\sigma_\delta$ belongs to $S_{1,\delta}^m$,
but $\sigma_\delta(X,D)$ is not bounded on $M^{p,q}(\R^n)$.

We next consider the case $q \le 2$.
In this case,
$m>-(1/q-1/2)\delta n$.
Since $q' \ge 2$ and $m>(1/q'-1/2)\delta n$,
by Theorem \ref{2.1} with $q \ge 2$,
we see that $\sigma_\delta(X,D)$ is not bounded on $M^{p',q'}(\R^n)$.
By duality and \eqref{(2.2)},
if $\sigma_\delta^*(X,D)$ is bounded on $M^{p,q}(\R^n)$,
then $\sigma_\delta(X,D)$ is bounded on $M^{p',q'}(\R^n)$.
Hence,
$\sigma_\delta^*$ belongs to $S_{1,\delta}^m$,
but $\sigma_\delta^*(X,D)$ is not bounded on $M^{p,q}(\R^n)$.
The proof is complete.



\begin{thebibliography}{20}
\bibitem{C-V}
A.P. Calder\'on and R. Vaillancourt,
{A class of bounded pseudo-differential operators},
Proc. Nat. Acad. Sci. U.S.A. 69 (1972), 1185-1187.


\bibitem{Feichtinger1}
H.G. Feichtinger,
{Banach spaces of distributions of Wiener's type and interpolation},
in: P. Butzer, B.Sz. Nagy and E. G\"orlich (Eds.),
Proc. Conf. Oberwolfach,
Functional Analysis and Approximation,
August 1980,
Int. Ser. Num. Math.,
Vol. 69,
Birkh\"auser-Verlag, Basel, Boston, Stuttgart, 1981, 153-165.

\bibitem{Feichtinger2}
H.G. Feichtinger,
{Modulation spaces on locally compact abelian groups},
Technical Report,
University of Vienna, Vienna, 1983.

\bibitem{Feichtinger3}
H.G. Feichtinger,
{Modulation spaces: Looking back and ahead},
Sampl. Theory Signal Image Process. 5 (2006), 109-140.


\bibitem{Grochenig}
K. Gr\"ochenig,
{Foundations of Time-Frequency Analysis},
Birkh\"auser, Boston, 2001.

\bibitem{Grochenig2}
K. Gr\"ochenig,
{Time-Frequency analysis of Sj\"ostrand's class},
to appear in Rev. Mat. Iberoam.

\bibitem{Grochenig-Heil}
K. Gr\"ochenig and C. Heil,
{Modulation spaces and pseudodifferential operators},
Integral Equations Operator Theory 34 (1999), 439-457.

\bibitem{Kumano-go}
H. Kumano-go,
{Pseudo-Differential Operators},
MIT Press, Cambridge, 1981.

%

\bibitem{Sjostrand}
J. Sj\"ostrand,
{An algebra of pseudodifferential operators},
Math. Res. L. 1 (1994), 185-192.


\bibitem{Tachizawa}
K. Tachizawa,
{The boundedness of pseudodifferential operators
on modulation spaces},
Math. Nachr. 168 (1994), 263-277.

\bibitem{Toft}
J. Toft,
{Continuity properties for modulation spaces,
with applications to pseudo-differential calculus-I},
J. Funct. Anal. 207 (2004), 399-429.


\bibitem{Triebel}
H. Triebel,
{Modulation spaces on the Euclidean $n$-spaces},
Z. Anal. Anwendungen 2 (1983), 443-457.
\end{thebibliography}
\end{document}